\documentclass[11pt,onecolumn]{article}

\usepackage{amsmath}
\usepackage{theorem}
\usepackage{amsfonts}
\usepackage{graphicx}
\usepackage{psfrag}
\usepackage{xcolor}
\usepackage{cancel} 

\renewcommand{\epsilon}{\varepsilon}

\newcommand{\ital}{\emph}

\newcommand{\KL}{\mathcal K\mathcal L}
\newcommand {\R}{\mathbb R}
\renewcommand{\P}{\mathcal P}

\newcommand{\defby}{:=}

\newcommand{\foral}{\forall\,}

\newcommand{\pp}{p\in\P}

\usepackage[normalem]{ulem} 

{\theorembodyfont{\rmfamily} }

\begin{document}

\title{\bf Commutation relations and stability of switched systems:\\ a personal history}
\author{Daniel Liberzon\\University of Illinois Urbana-Champaign}
\date{}
\maketitle

\begin{abstract}
This expository article presents an overview of research, conducted mostly between the mid-1990s and late 2000s, that explores a link between commutation relations among
a family of asymptotically stable vector fields and stability properties of the
switched system that these vector fields generate. This topic is viewed through the lens of the author's own involvement with it, by interspersing explanations of technical developments with personal reminiscences and anecdotes.
\end{abstract}

\begin{center}
	The author dedicates this article to himself on the occasion of his 50th birthday.
\end{center}

\section{How it all began for me}
{In January 1998, having just defended my PhD thesis, I came to Yale to do a postdoc with Steve Morse (for what eventually turned out to be 2.5 very enjoyable years). At that time Steve was getting interested in stability of switched systems, a subject I knew nothing about. A particular result that caught his attention was by Leonid Gurvits who, I believe, had given a seminar at Yale shortly before my arrival. This result---which we will examine in detail in Section~\ref{s-gurvits} below---established a link between stability of a switched linear system and the property that the Lie algebra generated by the individual system matrices is nilpotent. Steve knew that I was trained as a mathematician and that I did my PhD under Roger Brockett, who was a pioneer in the use of Lie algebras in control theory (see, e.g., \cite{brockett-sicon-1972,brockett-lie-algebras-73}). Grossly overestimating my knowledge of this subject, Steve suggested that I try to improve and generalize Gurvits' result.}

To be able to explain what happened next, I first need to more formally introduce switched systems and their stability properties. In what follows, I assume that the reader is familiar with differential equations and understands what the notation $\dot x=f(x)$ means, but not much else.

\section{Switched systems}
Suppose we are given a collection $f_p:\R^n\to\R^n$, $p\in\P$ of vector fields or, what is the same, a collection of dynamical systems (which we will call \ital{modes})
\begin{equation}\label{e-modes}
\dot x=f_p(x),\qquad p\in\P
\end{equation}
with state $x\in\R^n$. Here $\P$ is an index set, which we assume for simplicity to be finite throughout this article by setting $\P=\{1,\dots,m\}$ for some positive integer $m$, although this assumption is not really necessary for many of the results that we will discuss.

A second ingredient needed to define a switched system is a function of time that specifies which of the above modes is to be activated when. We denote this function by $\sigma:[0,\infty)\to\P$ and call it a \ital{switching signal}. While the total number of switches (discontinuities of $\sigma$) can be infinite, we ask that on every \ital{bounded} time interval the number of switches be finite (to avoid some unpleasant technical issues). Since $\P$ is a finite set, it is clear that $\sigma(\cdot)$ is a piecewise constant function. Therefore, we are dealing with time variation of a fairly restricted type.

The resulting dynamics can be written as
$$
\dot x(t)=f_{\sigma(t)}(x(t))
$$
or, more compactly, as
\begin{equation}\label{e-switched-nonlinear}
\dot x=f_\sigma(x).
\end{equation}
We call the system~\eqref{e-switched-nonlinear} a \ital{switched system}, and it will be our main object of study. While in writing~\eqref{e-switched-nonlinear} I am using the notation that Steve Morse adopted and showed to me in 1998, switched systems had been studied in the control theory literature for some years prior to that and, in view of their proximity to differential inclusions~\cite{aubin-cellina-inclusions-book} and discontinuous systems~\cite{filippov-book}, their history goes back several decades.

\section{Stability notions}\label{s-stability-notions}

Let us suppose that all the individual modes~\eqref{e-modes} share an equilibrium at the origin, i.e., $f_p(0)=0$ for all $p\in\P$, and that we are interested in (asymptotic) stability properties of this equilibrium for the switched system~\eqref{e-switched-nonlinear}. More precisely, we want to know when such stability properties hold uniformly over all (piecewise constant) switching signals $\sigma(\cdot)$. The simplest such property to define is \ital{global uniform exponential stability (GUES)}, which asks that all solutions of~\eqref{e-switched-nonlinear} satisfy
\begin{equation}
|x(t)|\le c e^{-\lambda t}|x(0)| \qquad \foral t\ge 0, \ \foral\sigma(\cdot)
\tag*{(GUES)}
\end{equation}
for some fixed positive constants $c$ and $\lambda$. (Here $|\cdot|$ is some chosen norm in $\R^n$, usually the Euclidean norm.)

Note that the previous inequality is of the form
\begin{equation}
|x(t)|\le \beta(|x(0)|,t) \qquad \foral t\ge 0, \ \foral\sigma(\cdot)
\tag*{(GUAS)}
\end{equation}
where $\beta(r,t):=ce^{-\lambda t}r$. The reason we labeled this more general property as ``GUAS" is that, with $\beta:[0,\infty)\times[0,\infty)\to [0,\infty)$ an arbitrary continuous function such that $\beta(\cdot,t)$ is 0 at 0 and increasing for each fixed $t$ and $\beta(r,\cdot)$ is decreasing to 0 for each fixed $r$, this gives us precisely \ital{global uniform asymptotic stability}, which is GUES but without the requirement that the solution bound decay exponentially in time or grow linearly in the norm of the initial state. Functions $\beta$ with the above properties are known as ``class $\KL$ functions"; for example, $\beta(r,t):=r^2/(1+t)$ is one such function. For the case of a single (non-switched) system, the above properties~(GUES) and~(GUAS) reduce to the classical global exponential stability and global asymptotic stability, respectively.

If GUAS still seems like too strong a requirement, we can consider its local version (by requiring the bound to hold only for initial states in some neighborhood of 0). We can also drop the uniformity over $\sigma$, and ask instead that for each \ital{fixed} switching signal $\sigma(\cdot)$, the corresponding (time-varying) system be asymptotically stable (globally or locally), or just ask that its solutions converge to 0 (and ignore Lyapunov stability). But we should admit that even this latter property---\ital{attractivity for each fixed $\sigma$}---is still very strong. Indeed, it is easy to construct examples where all modes are asymptotically stable but switching between them can lead to instability (see, e.g., \cite{survey}).

It may thus seem more reasonable---and more practically relevant---to ask for asymptotic stability to be preserved only for \ital{some}, but not all, switching signals, and to investigate what these switching signals are. So, why insist on asymptotic stability for \ital{all} switching signals? I have at least two answers to this question. One is that examining stability under arbitrary switching, and understanding possible instability mechanisms, paves the way to studying the more realistic scenario of stability under constrained switching. Another answer is that mathematically speaking, the problem of asymptotic stability under arbitrary switching turns out to be interesting and admits elegant solutions, as I try to demonstrate below.

It should be clear that a prerequisite for stability (in the asymptotic, exponential, or any other sense) is that each individual mode possess this stability property, since constant switching signals $\sigma\equiv p\in\P$ are allowed. This necessary condition will be assumed throughout.

\section{Common Lyapunov functions}\label{s-clf}

Readers familiar with basic stability theory for nonlinear systems can skim or skip the following two paragraphs, while others can consult a textbook such as~\cite{khalil-book-3ed} for further details.
Let us go back for a moment to a single (non-switched) system
\begin{equation}\label{e-nonswitched}
\dot x=f(x)
\end{equation}
with $x\in\R^n$ and $f(0)=0$. A classical approach known as \ital{Lyapunov's direct method} analyzes stability of~\eqref{e-nonswitched} with the help of a \ital{Lyapunov function}. As a candidate Lyapunov function, one takes a continuously differentiable function $V:\R^n\to \R$ such that $V(0)=0$ and $V(x)>0$ for all $x\ne 0$ (such $V$ are called \ital{positive definite}). Then one defines a new function, $\dot V(x):=({\partial V}/{\partial x})\cdot f(x)$. Its significance lies in the fact that the derivative of $V(x(t))$ along solutions of the system~\eqref{e-nonswitched} is given by $\frac d{dt}V(x(t))=\dot V(x(t))$. The key ideas here are that the function $\dot V$ can be obtained without the knowledge of the system's solutions, and that this allows us to reduce the analysis of the complex behavior of $x(t)\in\R^n$ to that of the scalar-valued quantity $V(x(t))$. Specifically, Lyapunov's stability theorem asserts that the system~\eqref{e-nonswitched} is stable (in the sense of Lyapunov) if $\dot V(x)\le 0$ for all $x$; asymptotically stable if $\dot V(x)< 0$ for all $x\ne 0$; and globally asymptotically stable if $V$ also has the property that $V(x)\to\infty$ whenever $|x|\to \infty$ (called \ital{radial unboundedness}). Exponential stability can be concluded under some additional structure, for example, if both $V$ and $\dot V$ are quadratic (or are sandwiched between quadratic bounds).

It is important to know that \ital{converse} Lyapunov theorems also exist, whereby stability of the system (asymptotic or exponential, local or global) can be used to prove the existence of a Lyapunov function with requisite properties. The associated formulas for the Lyapunov functions involve the system's solutions, and are thus not constructive. Let us denote by $\phi(t,x)$ the solution of~\eqref{e-nonswitched} at time $t$ corresponding to the initial condition $x(0)=x$. Under the assumption that~\eqref{e-nonswitched} is asymptotically stable, a construction due to Massera produces a Lyapunov function of the form
\begin{equation}\label{e-massera}
V(x)=\int_0^\infty G(|\phi(t,x)|)dt
\end{equation}
where $G:[0,\infty)\to[0,\infty)$ is a function with certain properties. In the particular case when~\eqref{e-nonswitched} is exponentially stable, we can take $G(r)=r^2$, resulting in
\begin{equation}\label{e-massera-es}
V(x)=\int_0^\infty |\phi(t,x)|^2 dt.
\end{equation}
Another construction, due to Kurzweil, does not involve integration and instead works with the function $g(x):=\inf_{t\le 0}|\phi(t,x)|$ and then defines
\begin{equation}\label{e-kurzweil}
V(x):=\sup_{t\ge 0} g(\phi(t,x))k(t)
\end{equation}
for an appropriate auxiliary function $k(\cdot)$. Kurzweil's construction works globally if asymptotic stability is global, while Massera's construction only works on a bounded region. On the other hand, Massera's construction is capable of providing an upper bound on the norm of the gradient ${\partial V}/{\partial x}$, which is important for perturbation analysis.

Let us return to the switched system~\eqref{e-switched-nonlinear}. Just like in the non-switched case, we can analyze its stability with the help of a Lyapunov function $V:\R^n\to \R$. But now we want $V$ to decay along solutions of each mode, i.e., we want to have $({\partial V}/{\partial x})\cdot f_p(x)<0$ for all $x\ne 0$ and all $p\in\P$. Moreover, we want a uniform lower bound on this decay rate over all modes. For a finite collection of modes (which is the situation considered here) the existence of such a mode-independent decay rate is automatic, but it still helps to write this property explicitly as follows:
\begin{equation}\label{e-common-Lyapunov}
\frac{\partial V}{\partial x}f_p(x)\le -W(x)<0 \qquad\foral x\ne 0,\ \foral p\in\P
\end{equation}
for some continuous function $W$. A function $V$ satisfying~\eqref{e-common-Lyapunov} is referred to as a \ital{common Lyapunov function} for our family of systems~\eqref{e-modes} or, with a slight abuse of terminology, for the switched system~\eqref{e-switched-nonlinear}. We see from~\eqref{e-common-Lyapunov} that when we start switching, $V(x(t))$ will always decay along solutions of~\eqref{e-switched-nonlinear} at the rate of at least $W(x(t))$. This is all that is needed to show uniform asymptotic stability, by arguing in the same way as in the proof of Lyapunov's classical stability theorem (the fact that $V(x(\cdot))$ is not differentiable at the switching times does not really affect the proof). If $V$ is also radially unbounded, we conclude GUAS. Converse Lyapunov theorems are also available, stating that the existence of a common Lyapunov function is necessary for GUAS. (Some of these results were actually derived in the more general setting of differential inclusions;
see, e.g., \cite[p.\ 188]{book} for further discussion and references.)

When trying to prove GUAS (or one of its variants), we typically have two choices: one is to directly analyze the trajectories in the time domain, and the other is to look for a common Lyapunov function. We will soon see examples of successful application of both of these approaches. However,
proving GUAS by either method
is very challenging in general,
and so we need to identify some additional system structure that can help us.

\section{Switched linear systems and their stability}\label{s-linear-stab}

Students of systems and control theory (particularly those in engineering disciplines rather than in mathematics) are typically first introduced to linear systems before being exposed to more general nonlinear systems. What makes the analysis of linear systems more manageable is the availability of an explicit formula for system solutions in terms of matrix exponentials, as well as computationally tractable stability conditions in terms of quadratic Lyapunov functions. Accordingly, a very widely studied special case of the switched system~\eqref{e-switched-nonlinear} is obtained by taking the individual modes to be linear, of the form $\dot x=A_px$ with each $A_p$ a real-valued $n\times n$ matrix. This leads to the \ital{switched linear system}
\begin{equation}\label{e-switched-linear}
\dot x=A_\sigma x.
\end{equation}
For reasons explained previously, we take the matrices $A_p$, $p\in\P$ to be Hurwitz (i.e., with eigenvalues having negative real parts).

For our purposes, there is one more feature that makes studying switched linear systems convenient. As shown by David Angeli around 1999~\cite{angeli-switched-note-scl}, the different stability notions that we introduced in Section~\ref{s-stability-notions} for the switched nonlinear system~\eqref{e-switched-nonlinear} actually all turn out to be equivalent in the linear case. Namely, for the switched linear system~\eqref{e-switched-linear} we have
\begin{equation}\label{e-equivalences}
\mbox{GUES}\ \Leftrightarrow \ \mbox{GUAS} \ \Leftrightarrow \ \mbox{(local) UAS}
\ \Leftrightarrow \ \mbox{attractivity for each } \sigma(\cdot).
\end{equation}
All the ``$\Rightarrow$" implications are of course obvious from the definitions (and do not require linearity). The converse implications are deduced with the help of the fact that solutions of~\eqref{e-switched-linear} scale linearly with initial conditions.\footnote{It is enough to know that the dependence on initial conditions is homogeneous of degree 1.} The middle ``$\Leftarrow$" implication follows immediately from this observation. The leftmost  ``$\Leftarrow$" implication is also not difficult to justify by noting that, if $\beta$ is the function on the right-hand side of~(GUAS) and if $T>0$ is such that, say, $\beta(1,T)\le 1/2$, then by linearity we must have $|x(T)|\le |x(0)|/2$,  $|x(2T)|\le |x(0)|/4$, and so on, hence the convergence is exponential. The remaining (rightmost)  ``$\Leftarrow$" implication is a deeper result, and relies on a general theorem about uniform attractivity proved by Sontag and Wang in~\cite{sontag-wang-tac-96}.

{I learned about the above equivalences from David Angeli when I met him during one of my visits to Rutgers University hosted by Eduardo Sontag, who during my postdoc years at Yale regularly invited me and other young researchers to come and interact. A few years later, when I published my switched systems book~\cite{book}, David's paper was still cited there as a manuscript submitted for publication. Apparently, the reviewers felt that the result was too simple and should be well known, probably buried somewhere in the Russian literature (although neither David nor the reviewers were able to locate a precise reference). Eventually, this result was added to the paper~\cite{angeli-leenheer-sontag-09} and this way it was finally published.}

We know from Section~\ref{s-clf} that the equivalent stability notions listed in~\eqref{e-equivalences} are also equivalent to the existence of a common Lyapunov function $V$ satisfying~\eqref{e-common-Lyapunov}. In general, constructing a common Lyapunov function is hard and there is no systematic procedure for doing it. But maybe the linear structure can help us? After all, for a single linear system $\dot x=Ax$, with $A$ a Hurwitz matrix, finding a Lyapunov function is straightforward. As is well known, for every matrix $Q=Q^T>0$ there is a unique matrix $P=P^T>0$ solving the Lyapunov equation
\begin{equation}\label{e-Lyap-eq}
PA+A^TP=-Q
\end{equation}
and this yields a quadratic Lyapunov function
\begin{equation}\label{e-qclf}
V(x)=x^TPx
\end{equation}
whose derivative along solutions is $-x^TQx$. For future reference, we note that $P$ is given by the explicit formula
\begin{equation}\label{e-Lyapunov-solution-explicit}
P=\int_0^\infty e^{A^Tt}Qe^{At}dt
\end{equation}
although in practice one just solves the linear system of equations~\eqref{e-Lyap-eq} directly. Remembering that the flow of $\dot x=Ax$ is given by $\phi(t,x)=e^{At}x$, we also see that for $Q=I$ this Lyapunov function construction
is a special instance of the one given by the formula~\eqref{e-massera-es}.

Now, for the switched linear system~\eqref{e-switched-linear} it seems reasonable to ask whether we can always search for a common Lyapunov function within the class of positive definite quadratic forms~\eqref{e-qclf}. This means finding a matrix $P=P^T>0$ that satisfies
\begin{equation}\label{e-LMIs}
PA_p+A_p^TP<0 \qquad \foral p\in\P
\end{equation}
or, more explicitly (but equivalently in view of the finiteness of $\P$),
\begin{equation}\label{e-qclf-PQ}
PA_p+A_p^TP\le -Q<0 \qquad \foral p\in\P.
\end{equation}
What is attractive about this possibility is that~\eqref{e-LMIs} is a system of \ital{linear matrix inequalities} (LMIs), and there are efficient numerical methods for solving them (see~\cite{lmi-book} for a comprehensive introduction, and~\cite{gradient-tac} for an alternative approach that handles the inequalities sequentially rather than simultaneously).

\begin{figure}[htbp]
\psfrag{a1}[][]{{$A_1 x$}} \psfrag{a2}[][]{{$A_2 x$}}
\centerline{\small \includegraphics[height=2in]{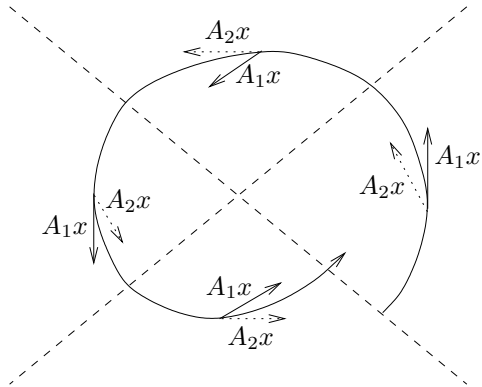}}
  \caption{\small Worst-case switching
\label{f-dayamartin}}
\end{figure}

Unfortunately, it turns out that working with quadratic common Lyapunov functions is not sufficient even for switched linear systems. In their 1999 paper~\cite{dayawansa-martin-commonlf-99}, Dayawansa and Martin gave an example of a switched linear system that is GUES but does not possess a quadratic common Lyapunov function.\footnote{A much older paper by Pyatnitskiy~\cite{pyatnitskiy-switched-1971} contains a different example that illustrates the same point, although it takes more effort to extract this observation from~\cite{pyatnitskiy-switched-1971} because that paper studies absolute stability of a class of time-varying systems and does not explicitly treat switched systems.} Their example involves constructing two $2\times 2$ Hurwitz matrices $A_1$ and $A_2$ for which the pair of inequalities~\eqref{e-LMIs} can be directly shown to be infeasible. The proof that the corresponding planar switched linear system is nevertheless GUES is more conceptually interesting, as it relies on an idea that will resurface prominently later in this article. This idea is to analyze the behavior
of the system
under ``worst-case switching,"
which in this particular case is described as follows. The linear vector fields $A_1 x$ and $A_2 x$ are
collinear on two lines passing through the origin (the dashed lines
in Figure~\ref{f-dayamartin}). In each conic region between these two lines, one
of the two vector fields points outwards relative to the other. The
worst-case switching strategy thus consists of following the vector
field that points outwards, with switches occurring on the two
lines. It can be verified that this produces a trajectory converging to
the origin, because the distance from the origin after one
rotation decreases (see Figure~\ref{f-dayamartin}). It is easy to see that all other
trajectories of the switched system must then also converge
to the origin. To conclude GUES, we can use the worst-case trajectory  to obtain
a uniform upper bound on the norm of all other trajectories (with the same initial state), or alternatively we can appeal to~\eqref{e-equivalences}.

The fact that GUES of a switched linear system is not always certifiable by a quadratic common Lyapunov function does not mean, of course, that such functions cannot still be useful in specific scenarios. In fact, we will see that they are quite useful when the Hurwitz matrices $A_p$, $p\in\P$ satisfy certain commutation relations.

\section{Commuting matrices}\label{s-commuting}

We are now finally ready to begin talking about what these commutation relations are and why they are relevant to the problem of stability under arbitrary switching. Our starting point is the following simple observation:

\smallskip
\emph{If the Hurwitz matrices $A_p$, $p\in\P$ commute pairwise, then the switched linear system~\eqref{e-switched-linear} is GUES.}
\smallskip

\noindent To see why this is true, take for simplicity $\P=\{1,2\}$. The commutativity condition of course just means that $A_1A_2=A_2A_1$, which we can also write as $[A_1,A_2]=0$ with the \ital{commutator}, or \ital{Lie bracket}, of two matrices defined as
\begin{equation}\label{e-Lie-bracket-matrices-def}
[A_1,A_2]:=A_1A_2-A_2A_1.
\end{equation}
For a switching signal that takes the value, say, 1 for $t_1$ units of time, then the value 2 for $s_1$ units of time, then 1 for $t_2$ units of time, and so on, the corresponding solution of~\eqref{e-switched-linear} at time $t$ is of the form
\begin{equation}\label{e-exps-before-rearranged}
x(t)= e^{A_2s_k} e^{A_1t_k}  \cdots e^{A_2s_2} e^{A_1t_2}e^{A_2s_1} e^{A_1t_1}x(0).
\end{equation}
Since $A_1$ and $A_2$ commute, their matrix exponentials commute as well. Hence, we can rearrange the above expression as\footnote{Here we are also using the fact that, e.g.,  $e^{A_1t_1}e^{A_1t_2}=e^{A_1(t_1+t_2)}$.}
\begin{equation}\label{e-exps-after-rearranged}
x(t)=e^{A_2(s_k+\cdots+s_2+s_1)}e^{A_1(t_k+\cdots+t_2+t_1)}x(0).
\end{equation}
As $t\to\infty$, at least one of the two sums (which represent the total activation times of modes 1 and 2) must also converge to $\infty$. Both $A_1$ and $A_2$ being Hurwitz, this ensures that at least one of the matrix exponentials in~\eqref{e-exps-after-rearranged} converges to 0. Therefore, the switched system is attractive for every switching signal, and we know from~\eqref{e-equivalences} that this implies GUES. This result (more precisely, the attractivity part) was briefly noted in the paper~\cite{narendra-balakrishnan-commuting}, whose main contribution we are about to discuss, but it was almost certainly known much earlier.

The expressions~\eqref{e-exps-before-rearranged} and~\eqref{e-exps-after-rearranged} are equivalent because the flows of the two linear systems $\dot x=A_1x$ and $\dot x=A_2 x$ commute. This equivalence can also be interpreted as follows: \ital{every state that can be reached from $x(0)$ with an arbitrary number of switches can also be reached (at the same time $t$) with at most one switch.} Thinking about the commutativity condition in terms of reachability with a bound on the number of switches suggests a path that will prove very fruitful for us later on, although initially it is easy to overlook it (as I certainly did when I first learned about this result).

As we discussed in Section~\ref{s-clf}, an alternative to a direct trajectory-based stability proof---such as the one just given---is to search for a common Lyapunov function. And when dealing with a switched linear system, we have a natural candidate in the form of a quadratic Lyapunov function~\eqref{e-qclf}, even though we cannot be sure a priori that a quadratic common Lyapunov function exists (see Section~\ref{s-linear-stab}). Fortunately, for the case of commuting matrices such a function does exist, and it can be constructed in an elegant way proposed by Narendra and Balakrishnan in~\cite{narendra-balakrishnan-commuting}. With $\P=\{1,\dots, m\}$ as before, one iteratively solves the sequence of Lyapunov equations
\begin{equation}\label{e-nb94}
\begin{split}
P_1A_1+A_1^TP_1&=-I,\\
P_2A_2+A_2^TP_2&=-P_1,\\
&\vdots\\
P_mA_m+A_m^TP_m&=-P_{m-1}
\end{split}
\end{equation}
(placing each $P_i$ obtained at step $i$ on the right-hand side of the equation to be solved at step $i+1$). One
then defines
\begin{equation}\label{e-VPm}
V(x):=x^TP_mx.
\end{equation}
It is obvious that this is a Lyapunov function for the $m$th mode, $\dot x=A_m x$, but commutativity implies that it is a Lyapunov function for all the other modes as well. One way to see this is to repeatedly apply the formula~\eqref{e-Lyapunov-solution-explicit} for the solution of the Lyapunov equation~\eqref{e-Lyap-eq} to express $P_m$ as a nested integral involving products of matrix exponentials:
$$
P_m= \int_0^\infty e^{A_m^Tt_m}\dots\Big(\int_0^\infty
e^{A_1^Tt_1} e^{A_1t_1}dt_1\Big)\dots e^{A_mt_m}dt_m.
$$
Since these matrix exponentials commute, we can arbitrarily reorder the integrals. This shows that the final matrix $P_m$ does not depend on the ordering of the modes in~\eqref{e-nb94}, giving the desired claim. It is also fairly easy (especially for the case $m=2$) to reach the same conclusion by manipulating the Lyapunov equations in~\eqref{e-nb94}, as done in~\cite{narendra-balakrishnan-commuting}.

{It is interesting to note that the work~\cite{narendra-balakrishnan-commuting} was done at Yale, by a different group and a few years before I arrived there to work with Steve Morse. I guess one could say that switched systems were in the air there during that time period.}

\section{Matrices generating a nilpotent Lie algebra}\label{s-gurvits}

We have now arrived at the point in our story where Leonid Gurvits comes in, with the result that I mentioned at the very beginning of the article. I first need to define the terms appearing in the title of this section. A (matrix) \ital{Lie algebra} generated by our collection of real-valued $n\times n$ matrices $A_p$, $p\in\P$ is the smallest set of matrices containing $\{A_p, p\in\P\}$ which is closed under addition and scalar multiplication (hence it is a vector space of matrices) and also under the Lie bracket operation $[\cdot,\cdot]$. In the commuting case considered in Section~\ref{s-commuting}, this Lie algebra is simply the linear span of the original matrices. In general, however, it must also contain their pairwise Lie brackets $[A_i,A_j]$, the iterated (second-order) Lie brackets $[A_i,[A_j,A_k]]$, and so on. Of course, even though the number of Lie brackets to account for is potentially infinite, the resulting vector space has dimension at most $n^2$, so the process of adding new Lie brackets linearly independent from the previous ones will terminate after finitely many steps.

The Lie algebra is \ital{nilpotent} if there is a positive integer $k$ such that all $k$th-order Lie brackets are 0. For example, for $m=2$ (two matrices) and $k=2$ (second-order nilpotency) this means that
\begin{equation}\label{e-2nd-order-nilpotent-linear}
[A_1,[A_1,A_2]]=[A_2,[A_1,A_2]]=0.
\end{equation}
This scenario represents the next natural step after the commuting case, in the sense that the first-order Lie bracket $[A_1,A_2]$ is the only one that needs to be added.

In the 1995 paper~\cite{gurvits-nilpotent-laa95}, Gurvits made the following conjecture:\footnote{Actually Gurvits worked in discrete time, so here we are paraphrasing his statements in the continuous-time setting. The difference disappears if we restrict the switching times to be integer multiples of a fixed positive number.}

\smallskip
\emph{If the Hurwitz matrices $A_p$, $p\in\P$ generate a nilpotent Lie algebra, then the switched linear system~\eqref{e-switched-linear} is GUES.}
\smallskip

\noindent Gurvits only proved this claim for the special case of two matrices generating a second-order nilpotent Lie algebra, represented by the condition~\eqref{e-2nd-order-nilpotent-linear}. His method was the following. We saw that for two commuting matrices, an arbitrary product of matrix exponentials as in~\eqref{e-exps-before-rearranged} representing a solution of the switched system at time $t$ can be equivalently written as a product of the form $e^{A_2\tau_2}e^{A_1\tau_1}$ with $\tau_1+\tau_2=t$, appearing in~\eqref{e-exps-after-rearranged}. For the second-order nilpotent case, Gurvits showed that an arbitrary product from~\eqref{e-exps-before-rearranged} can be similarly represented as a shorter product, but this time of \ital{five} rather than two exponentials:
\begin{equation}\label{e-gurvits-factorization}
x(t)=e^{A_2\tau_5}e^{A_1\tau_4}e^{A_2\tau_3}e^{A_1\tau_2}e^{A_2\tau_1}x(0).
\end{equation}
In the terminology of~\cite{gurvits-nilpotent-laa95}, the semigroup generated by the two matrices can be factored into a product of five semigroups, each generated by one of the matrices. (In the commuting case, we have a product of two semigroups.)
Here the numbers $\tau_1,\dots,\tau_5$ do not necessarily add up to $t$, but still have the property that at least one of them must go to $\infty$ as $t\to \infty$. We can thus conclude GUES in the same way as in the commuting case. We can also interpret the decomposition~\eqref{e-gurvits-factorization} in terms of reachability with a bound on the number of switches, similarly to how we did in the commuting case: \ital{every state that can be reached from $x(0)$ with an arbitrary number of switches can also be reached (not necessarily at the same time $t$) with at most {four} switches.}

To prove~\eqref{e-gurvits-factorization}, Gurvits used the \ital{Baker-Campbell-Hausdorff
formula}
$$
e^{A_1}e^{A_2}=e^{A_1+A_2+\frac 12[A_1,A_2]+\frac
1{12}([A_1,[A_1,A_2]]+[A_2,[A_1,A_2]])+\dots}
$$
which, in the case of~\eqref{e-2nd-order-nilpotent-linear}, simplifies to just
$
e^{A_1}e^{A_2}=e^{A_1+A_2+\frac 12[A_1,A_2]}.$
It follows that an arbitrary product of matrix exponentials as in~\eqref{e-exps-before-rearranged} can be written as $e^{m_1A_1+m_2A_2+\frac 12 m_3[A_1,A_2]}$ for suitable coefficients $m_1,m_2,m_3$.
Using this fact, Gurvits analyzed the effect of an ``elementary shuffling", i.e., of switching the order of two adjacent matrix exponentials in a product, and showed that products of five terms as in~\eqref{e-gurvits-factorization} are sufficient to generate all possible values of $m_1,m_2,m_3$.

{As for the general nilpotent case (Lie algebra nilpotent of order $k\ge 2$ generated by $m\ge 2$ matrices), Gurvits wrote that it ``looks quite
hopeful, but requires more sophisticated combinatorics". This is exactly what Steve Morse suggested I try to work out. Of course, it is not important that the factorization have 5 terms as in~\eqref{e-gurvits-factorization}; as long as we have some fixed bound $N$ on the number of terms in the product, the stability proof goes through.}

{Initially I was optimistic about finding such an $N$ by using manipulations similar to the ones employed by Gurvits.
After spending some time trying to do this myself, without success, I asked for help. One person to whom I described the problem was Gerry Schwarz, an expert in Lie theory at Brandeis University's math department, where I received my PhD. After a few weeks Gerry sent me an email, announcing that he was very close to a solution and asking if I was still interested in one. I responded with an emphatic ``yes" and then never heard from him again. Another mathematician I consulted was George Seligman at Yale, who was kind enough to meet with me several times and educate me about various aspects of Lie algebras. While he did not have a direct solution to my problem, he eventually helped me stumble upon a different approach which not only led to a positive result, but in fact applied to a larger class of systems than the one studied by Gurvits.}

\section{Matrices generating a solvable Lie algebra}\label{s-solvable}

A well-known class of Lie algebras that contains all nilpotent Lie algebras is that of \ital{solvable} Lie algebras. While nilpotent Lie algebras are characterized by all Lie brackets of sufficiently high order being 0, in a solvable Lie algebra only Lie brackets of sufficiently high order having a certain structure must be 0. To make this a bit more precise, recall that 2nd-order nilpotency means that Lie brackets of the form $[A_i,[A_j,A_k]]$ vanish; 3rd-order nilpotency means that Lie brackets of the form $[A_i,[A_j,[A_k, A_\ell]]]$ vanish; and so on. Here at step $k$ one considers Lie brackets of matrices obtained at step $k-1$ with matrices from the original Lie algebra. By contrast, when defining solvability, at step $k$ one only considers Lie brackets among the matrices from step $k-1$; for example, for $k=2$ one looks at Lie brackets of the form $[[A_i,A_j],[A_k,A_\ell]]$. The latter approach singles out a subset of the Lie brackets included when using the former approach. Therefore, every nilpotent Lie algebras is solvable, while some solvable Lie algebras are not nilpotent.
During one of my conversations with George Seligman, he called my attention to solvable Lie algebras and pointed out their classical characterization known as \ital{Lie's theorem}. This theorem says that matrices in a solvable Lie algebra can be simultaneously brought to an upper-triangular form by some linear (generally complex-valued) change of coordinates.

{At this point I must bring in another major character in this story: Jo\~ao Hespanha, who at that time was finishing up his PhD studies under Steve Morse. We overlapped at Yale for only one semester, but that was enough to establish collaboration on a range of topics which continues on and off to this day. Stability of switched systems was an area in which Jo\~ao had already been working with Steve for some time before I arrived, motivated primarily by problems in switching adaptive control. So, when I mentioned to Jo\~ao solvability and triangular structure, in the context of trying to generalize Gurvits' approach, this immediately rang a bell for him. He had previously encountered switching among stable linear systems with triangular matrices, and he knew that such switched linear systems are always stable.}

To see why this is true, let us consider the case when $\P=\{1,2\}$ and
$x\in\R^2$. Let the two matrices be
\begin{equation}\label{e-triangular-2example}
A_1:=\begin{pmatrix}
  -a_1 & b_1\\
  0 & -c_1
\end{pmatrix},\qquad
A_2:=\begin{pmatrix}
  -a_2 & b_2 \\
  0 & -c_2
\end{pmatrix}.
\end{equation}
Suppose for simplicity that their entries are real (the case of
complex entries requires some care but the extension is not
difficult). Since the eigenvalues of these matrices have negative
real parts, we have $a_i,c_i>0$, $i=1,2$. Now, consider the
switched linear system $\dot x=A_\sigma x$. The second component of $x$
satisfies the equation
$\dot x_2=-c_\sigma x_2.$ Therefore, $x_2$ decays to zero
exponentially fast for every $\sigma(\cdot)$, at the rate corresponding to
$\min\{c_1,c_2\}$. The first component of $x$ satisfies the
equation $\dot x_1=-a_\sigma x_1 + b_\sigma x_2.$ This can be
viewed as the exponentially stable system $\dot x_1=-a_\sigma x_1$
perturbed by the exponentially decaying input $b_\sigma x_2$. Thus
$x_1$ also converges to zero exponentially fast. It is not hard to
extend this argument to more than two matrices of arbitrary dimension, proceeding from the bottom component of $x$ upward. As before, GUES follows by virtue of~\eqref{e-equivalences}.

An alternative to the above direct stability proof, as we know, consists in constructing a common
Lyapunov function. It turns out that in
the present case of triangular matrices, it is possible to find a quadratic common
Lyapunov function of the form~\eqref{e-qclf}, with $P$ a diagonal
matrix. We illustrate this again on the example of the two
matrices~\eqref{e-triangular-2example}. Let us look for $P$ taking
the form
$$
P=\begin{pmatrix}
  d_1 & 0 \\
  0 & d_2
\end{pmatrix}
$$
where $d_1,d_2>0$. A straightforward calculation gives
$$
-A_i^TP-PA_i=\begin{pmatrix}
  2d_1a_i& -d_1b_i \\
  -d_1b_i & 2d_2c_i
\end{pmatrix},\qquad i=1,2.
$$
To ensure that this matrix is positive definite, we can first pick
an arbitrary $d_1>0$, and then choose $d_2>0$ large enough to have
$
4d_2d_1a_ic_i-d_1^2b_i^2>0$, $i=1,2
$.
Again, it is easy to see how this iterative construction can be extended to higher dimensions.\footnote{A slightly different approach is to rescale the basis vectors to make the off-diagonal elements (the $b_i$'s in the present example) as small as desired, so that $P$ can be taken to be the identity matrix (see, e.g., \cite[pp.\ 203--204]{arnold-ode-book}).}

We arrive at the following result:

\smallskip
\emph{If the Hurwitz matrices $A_p$, $p\in\P$ generate a solvable Lie algebra, then the switched linear system~\eqref{e-switched-linear} is GUES.}
\smallskip

\noindent {Of course, it is clear from the preceding discussion that this is just a straightforward combination of two ingredients. The first one is the classical Lie's theorem, found in any textbook on Lie algebras, which gives the triangular form. The second ingredient is the observation that triangular form guarantees stability under arbitrary switching. Although I was initially unaware of this second result, Jo\~ao knew about it and in fact it had been documented in the literature. In particular, the paper~\cite{cohen-lewkowicz-triangular-laa-97}, published one year before our investigation, and the paper~\cite{shorten-narendra-cdc98}, written at about the same time (and the same place) as ours, both mention essentially the same Lyapunov function construction as the one given above.}

{On the other hand, it can be argued that the above result is greater than the sum of its parts. Indeed, the Lie-algebraic condition can be checked by performing a finite number of computations with the original matrices. It tells us that there exists a basis in which these matrices take the triangular form, but we do not need to actually find such a basis. And compared with Gurvits' approach described in Section~\ref{s-gurvits}, the present argument is much simpler and gives a stronger claim.}

{We summarized the above findings in a paper that we submitted to Systems and Control Letters. At first it was rejected---the reviewers thought that the result was trivial. Although we did not entirely disagree, we insisted that the paper still had value, and eventually it was published~\cite{lie}. In spite of (or maybe thanks to?) its almost embarrassing simplicity, the paper became quite highly cited. Unbeknown to me at the time, there was another reason to feel less than proud about this paper, but more on that later.}

It was noted in~\cite{haimovich-braslavsky-tac11} that the Lie-algebraic stability condition discussed in this section can also be used for control design purposes. Namely, given a family of linear control systems $\dot x=A_px+B_pu$, $\pp$, one can search for feedback gains $K_p$, $\pp$ such that the closed-loop matrices $A_p+B_pK_p$ are Hurwitz and generate a solvable Lie algebra or, equivalently, are simultaneously triangularizable. The closed-loop switched linear system will then be GUES. The paper~\cite{haimovich-braslavsky-tac11} describes an algorithm for finding such stabilizing feedback gains.

\section{More general matrix Lie algebras}\label{s-compact-levi-factor}

{The above line of research was continued in my joint work with Andrei Agrachev. Before describing this work, let me briefly digress to explain the role that Andrei Agrachev has played in my academic life. Third-year undergraduate students of mathematics at Moscow University had to choose an area of specialization. The process involved us listening to presentations made by professors in the department about their research. One of them was Agrachev, who spoke about a geometric approach to nonlinear controllability. I was immediately captivated, and this led me to choose control theory as my specialization area and Agrachev as my undergraduate research advisor. After working with him for about a year I went to graduate school in the United States, where I was still studying control theory but the specific topics were different, and for a while Agrachev and I lost contact. Then in 1998, after I defended my thesis and started attending international control conferences, we met again and reconnected. Having just recently written the paper~\cite{lie}, I showed a preprint to Agrachev. The next day, he told me that he had read it and could see how to generalize it.}

{Once I was back at Yale, our collaboration proceeded mostly by Agrachev sending me very short emails and me spending hours in the Yale math library deciphering them. The results of this work are documented in our 2001 paper~\cite{agrachev}.} That paper makes use of deeper results from the theory of Lie algebras compared to~\cite{lie}, and so here I will limit myself to a brief informal summary. If the Lie algebra generated by the matrices $A_p$, $p\in\P$ is not solvable, it contains a maximal solvable subalgebra ${\mathfrak{r}}$, called the \ital{radical}.\footnote{More precisely, ${\mathfrak{r}}$ is a maximal solvable ideal.} Every matrix $A_p$ can then be written as a sum $A_p=R_p+S_p$, where $R_p\in {\mathfrak{r}}$ and $S_p$ lies in a complementary subalgebra, which we call ${\mathfrak{s}}$. Suppose that all matrices in ${\mathfrak{s}}$ have purely imaginary eigenvalues (i.e., they are essentially rotation matrices); the subalgebra ${\mathfrak{s}}$ is then said to be \ital{compact}. What we showed in~\cite{agrachev} is that under this condition, the matrices $S_p\in {\mathfrak{s}}$ do not affect stability of the switched system.
We can paraphrase the resulting stability criterion as follows:

\smallskip
\emph{If the Hurwitz matrices $A_p$, $p\in\P$ generate a ``solvable plus compact" Lie algebra, then the switched linear system~\eqref{e-switched-linear} is GUES.}
\smallskip

A quadratic common Lyapunov function also exists in this case, although the proof of this fact given in~\cite{agrachev} is not nearly as constructive as in the solvable case. Moreover, it turns out that the above sufficient
condition for stability is
the strongest one that can be obtained by working solely with the Lie algebra. Indeed, we proved in~\cite{agrachev} that if the Lie algebra is not ``solvable plus compact" then it
can always be generated by a family of Hurwitz
matrices (which might be different from $A_p$, $p\in\P$) such
that the corresponding
switched linear system is not stable. Thus we have in some sense reached the end of the road in formulating stability conditions for switched linear systems in terms of commutation relations between their matrices. Of course, it is still possible to obtain stronger results by bringing in other tools (see Section~\ref{s-robustness} below for some further discussion on this).

\section{Commuting nonlinear vector fields}\label{s-commuting-nonlinear}

Let us now go back to the switched nonlinear system~\eqref{e-switched-nonlinear} generated by the family of vector fields~\eqref{e-modes}, which we assume to share a globally asymptotically stable equilibrium at the origin. In light of the previous developments for switched linear systems, it is natural for us to first examine the situation where these vector fields commute. This is the same as saying that the corresponding flows $\phi_p(\cdot,x)$ commute, where $\phi_p(t,x)$ denotes the solution at time $t$ of the system $\dot x=f_p(x)$ with initial condition $x(0)=x$. For smooth vector fields, this property is captured by the fact that their \ital{Lie brackets} defined by
\begin{equation}\label{e-Lie-bracket-nonlinear-def}
[f_p,f_q](x):=\dfrac{\partial f_q(x)}{\partial x}f_p(x)-
\dfrac{\partial f_p(x)}{\partial x}f_q(x)
\end{equation}
equal 0 for all $p,q\in\P$. For linear vector
fields $f_p(x)=A_px$ the right-hand side becomes
$(A_qA_p-A_pA_q)x$, which is consistent with the definition of the
Lie bracket of two matrices except for the
difference in sign. The following is a generalization of the result from Section~\ref{s-commuting}:

\smallskip
\emph{If the globally asymptotically stable vector fields $f_p$, $p\in\P$ commute pairwise, then the switched nonlinear system~\eqref{e-switched-nonlinear} is GUAS.}
\smallskip

To see why this is true, we can try arguing as in the linear case. Take $\P=\{1,2\}$ for simplicity, and consider a switching signal that takes the value, say, 1 for $t_1$ units of time, then the value 2 for $s_1$ units of time, then 1 for $t_2$ units of time, and so on. Since the two flows commute, the solution of~\eqref{e-switched-nonlinear} at time $t$ can be written as
\begin{equation}\label{e-one-switch-comm-nonlinear}
x(t)= \phi_2(s_k+\cdots+s_2+s_1, \phi_1(t_k+\cdots+t_2+t_1,x(0)))
\end{equation}
where the two sums represent the total activation times of the two modes. As $t\to\infty$, at least one of these sums must converge to $\infty$. Since both vector fields are globally asymptotically stable, we conclude that $x(t)\to\infty$.

Note that this argument falls short of establishing GUAS, because we no longer have the equivalences~\eqref{e-equivalences}. We must show the existence of a class $\KL$ function $\beta$ as in~(GUAS). For this, we can use the fact that each globally asymptotically stable mode has its own class $\KL$ function $\beta_p$ supplying the upper bound $|\phi_p(t,x)|\le\beta_p(|x|,t)$. We also know from~\eqref{e-one-switch-comm-nonlinear} that if a state can be reached at time $t$ from $x(0)$ with an arbitrary number of switches, then it can be reached at the same time $t$ with at most one switch. Combining these two properties, we can write $|x(t)|\le \beta_p(\beta_q(|x(0)|,t-\tau),\tau)$ where $\tau$ is the time of the switch and $(p,q)$ is either $(1,2)$ or $(2,1)$. Of the two time arguments $\tau$ and $t-\tau$, one is at least $t/2$; replacing it by $t/2$ and replacing the other time argument by 0 can only increase the upper bound. Taking the maximum over the different cases mentioned, we can easily obtain a function $\beta$ certifying GUAS.
A general construction of such a class $\KL$ function (for $m\ge 2$ modes) was given by Mancilla-Aguilar in the 2000 paper~\cite{mancilla-commuting-tac}; it proceeds a bit differently, by relying on some known properties of class $\KL$ functions and using induction on $m$.

In Section~\ref{s-commuting} we saw how a quadratic common Lyapunov function for a family of commuting exponentially stable linear systems can be constructed by the iterative procedure~\eqref{e-nb94}--\eqref{e-VPm}. We also noted that the standard quadratic Lyapunov function for a single exponentially stable linear system, obtained by solving the Lyapunov equation via~\eqref{e-Lyap-eq}--\eqref{e-Lyapunov-solution-explicit}, can be viewed as a special case of the Lyapunov function~\eqref{e-massera-es} for an exponentially stable nonlinear system. It is then natural to try to generalize the iterative construction of a common Lyapunov function to the case of commuting exponentially stable nonlinear vector fields. This was done by Hyungbo Shim and coworkers in~\cite{shim-...-siamcca-98}. With $\P=\{1,\dots,m\}$ again, one iteratively defines the functions
\begin{equation}\label{e-hyungbo}
\begin{split}
  V_1(x)&\defby \int_0^T|\phi_1(t,x)|^2dt,
  \\
  V_i(x)&\defby \int_0^TV_{i-1}(\phi_i(t,x))dt,\qquad i=2,\dots,m
\end{split}
\end{equation}
for a sufficiently large $T\le\infty$.  Then $V_m$
is a common Lyapunov function, at least locally; it is a global common Lyapunov function when the functions $f_p$, $p\in\P$ are
globally Lipschitz. For the case of linear systems $f_p(x)=A_px$ we
recover the construction of Section~\ref{s-commuting} upon setting $T=\infty$.

{I learned about the procedure~\eqref{e-hyungbo} from a poster that Hyungbo Shim presented at the 1998 SIAM Conference on Control and its Applications, which incidentally was the first international conference for both of us. After that initial meeting in front of Hyungbo's poster, a full decade would pass until he and I started collaborating (on topics not directly related to the present discussion).}

Note that an alternative to Hyungbo's approach is to employ Lyapunov's indirect (first) method. Namely,
consider the Jacobian matrices
\begin{equation}\label{e-jacobians}
A_p:=\frac{\partial f_p}{\partial x}(0),\qquad \pp.
\end{equation}
Lyapunov's indirect method tells us that the matrices $A_p$ are Hurwitz if (and only
if) the vector fields $f_p$ are (locally) exponentially stable.
Moreover, it can be shown that if the vector fields $f_p$ commute, then the
matrices $A_p$ also  commute.
(The converse does not necessarily hold, so the latter commutativity condition is actually weaker.) A
quadratic common Lyapunov function for the linearized systems $\dot x=A_px$ can thus be
constructed as explained in Section~\ref{s-commuting}, and it serves as a common
Lyapunov function for the original family of commuting nonlinear
vector fields. It is, however, only a local common Lyapunov function, and it only guarantees local uniform exponential stability of the switched nonlinear system~\eqref{e-switched-nonlinear}.

{In 2002 I was teaching a graduate course on switched systems that I had recently introduced at the University of Illinois. One of the students in the class was Linh Vu, who had just started his graduate studies with me. The course included final projects, which typically consisted in reading and presenting research articles, but Linh actually made an original research contribution. After I described Hyungbo Shim's Lyapunov function construction~\eqref{e-hyungbo} in class, Linh decided to develop analogous constructions for asymptotically but not necessarily exponentially stable commuting nonlinear systems. He proposed two approaches based on the classical Lyapunov function constructions of Massera and Kurzweil discussed in Section~\ref{s-clf}.} The first construction is based on~\eqref{e-massera} takes the form
\begin{align*}
  V_1(x)&\defby \int_0^\infty G(|\phi_1(t,x)|)dt,
  \\
  V_i(x)&\defby \int_0^\infty V_{i-1}(\phi_i(t,x))dt,\qquad i=2,\dots,m
\end{align*}
with the function $G$ coming from a nontrivial multivariable extension of a lemma due to Massera.\footnote{Interestingly, the Wikipedia entry for ``Massera's lemma" cites Linh's result.} The second construction is based on~\eqref{e-kurzweil} and takes the form
\begin{align*}
  V_1(x)&\defby \sup_{t\ge 0} g(\phi_1(t,x))k(t),
  \\
  V_i(x)&\defby \sup_{t\ge 0} V_{i-1}(\phi_i(t,x))k(t),\qquad i=2,\dots,m
\end{align*}
where $g(x)$ is the infimum norm of the backward-in-time solutions from $x$ over all switching signals, and $k(\cdot)$ is a suitable function.
These results became Linh's MS thesis and are documented in the paper~\cite{linh-commuting}, which also contains a readable and self-contained account of the relevant background material.

\section{Beyond commuting nonlinear vector fields: early attempts}

In view of the preceding developments, the next logical case to consider is when
the globally asymptotically stable vector fields $f_p$, $p\in\P$ generate a nilpotent or solvable Lie algebra. Here, the notions of a Lie algebra and its nilpotency and solvability are defined along the same lines as in Sections~\ref{s-gurvits} and~\ref{s-solvable}, except that instead of $n\times n$ matrices with the Lie bracket defined by~\eqref{e-Lie-bracket-matrices-def} we work with vector fields on $\R^n$ and the Lie bracket defined by~\eqref{e-Lie-bracket-nonlinear-def}.

We could, alternatively, inspect the Lie algebra generated by the Jacobian matrices~\eqref{e-jacobians}. If these matrices are Hurwitz and if the Lie algebra is solvable, then we know from Section~\ref{s-solvable} that a quadratic common Lyapunov function exists for the linearized systems $\dot x=A_px$, and then exactly as in Section~\ref{s-commuting-nonlinear} we can conclude local uniform exponential stability of the switched nonlinear system~\eqref{e-switched-nonlinear}. However, for the matrices $A_p$ to be Hurwitz, we need to assume that the vector fields $f_p$ are exponentially stable (at least locally). And overall this approach is not very interesting, as it is just an application of Lyapunov's indirect method. It seems much more intriguing to ask how the structure of the Lie algebra generated by
the original nonlinear vector fields $f_p$, $p\in P$ is related to (potentially global)
stability properties of the switched
nonlinear system~\eqref{e-switched-nonlinear}.

We saw in Section~\ref{s-solvable} that in the linear setting, solvability of the matrix Lie algebra implies simultaneous triangularizability (Lie's theorem), which in turn directly leads to uniform exponential stability of the switched linear system. It is tempting to try to carry out a similar program for the nonlinear vector fields $f_p$, $p\in\P$, where by the upper-triangular structure we now mean that they take the form
$$
f_p(x)=
\begin{pmatrix}
  f_{p1}(x_1,x_2,\dots,x_n)\\
  f_{p2}(x_2,\dots,x_n)\\
  \vdots\\
  f_{pn}(x_n)
\end{pmatrix}.
$$
Initially I was encouraged by the fact that there do exist nonlinear versions of Lie's theorem, which provide Lie-algebraic
conditions under which a family of nonlinear systems can be
simultaneously triangularized~\cite{crouch-triangular-siam-81,kawski-triangular-nilpotent,marigo-triangular-solvable-cdc99}.
These results unfortunately rely on some technical assumptions
that do not hold in our context, but let us ignore such details here. A more serious difficulty arises, however, when we try to explore the triangular structure for the purpose of establishing stability. When proving stability of the switched linear system generated by the matrices~\eqref{e-triangular-2example}, we used the fact that the state of an exponentially stable linear system perturbed by an input converging to 0 must converge to 0. For nonlinear systems, this ``converging-input-converging-state" property---which is a consequence of the well-known \ital{input-to-state stability (ISS)} property introduced by Sontag in~\cite{sontag-coprime}---is known not to be true in general. As an example, consider the system
\begin{equation}\label{e-alpha=1/2}
 \begin{split}
  \dot x_1&=-x_1+x_1^2x_2,\\
  \dot x_2&=-x_2.
\end{split}
\end{equation}
Even though $x_2\to 0$ exponentially fast, for sufficiently large initial conditions $x_1$ escapes to infinity in finite time (see, e.g., \cite[p.\ 8]{kkk-book} or \cite[p.\ 44]{book} for details).

{As I already mentioned in Section~\ref{s-linear-stab}, when I was doing a postdoc at Yale and working on this problem I met David Angeli who was interested in switched systems as well. During one of our discussions, I mentioned to David that I was trying to see if triangular structure might be helpful for showing GUAS of a switched nonlinear system. Any remaining hope I might have still had was promptly shattered by the following nice counterexample that David suggested.}
Let $\P=\{1,2\}$, and consider the upper-triangular vector fields
$$
f_1(x)=
\begin{pmatrix}
  -x_1+2\sin^2(x_1)x_1^2x_2\\
  -x_2
\end{pmatrix},\qquad f_2(x)=
\begin{pmatrix}
  -x_1+2\cos^2(x_1)x_1^2x_2\\
  -x_2
\end{pmatrix}.
$$
The systems $\dot x=f_1(x)$ and $\dot
x=f_2(x)$ are globally asymptotically stable, as one can easily verify by examining the behavior of their solutions. Nevertheless, the switched system $\dot x=f_\sigma(x)$ is not GUAS. Indeed, if it were GUAS, then we know from Section~\ref{s-clf} that $f_1$ and $f_2$ would share a common Lyapunov function. This Lyapunov function would in turn certify asymptotic stability of every ``convex combination" $\dot x=\alpha f_1(x)+(1-\alpha)f_2(x)$ of the two systems, where $\alpha\in[0,1]$. But for $\alpha=1/2$ we recover the unstable system~\eqref{e-alpha=1/2}, reaching a contradiction.

Our earlier discussion indicates that the problem here is that the $x_1$-dynamics are not ISS with respect to $x_2$ (viewed as an input). By imposing such ISS assumptions either on the switched system or on the individual modes, sufficient conditions for stability of switched triangular nonlinear systems can indeed be obtained, as David and I showed in~\cite{angeli}. (The above counterexample appears in the same paper.) However, such additional assumptions take us quite far from our original goal of formulating stability conditions in terms of the Lie algebra generated by the vector fields $f_p$, $\pp$.

We see that up to this point, all
attempts to formulate global asymptotic stability criteria valid
beyond the commuting nonlinear case were unsuccessful;
the methods employed to obtain the corresponding results for
switched linear systems do not apply, and an altogether different approach seems to
be required. I presented this as an open problem at a special session of the 2002 Mathematical Theory of Networks and Systems (MTNS) conference (it was included in the conference proceedings, and later published more formally as a book chapter~\cite{lie-open-problem}). By that time, I had mostly given up on this particular research direction and turned my attention to other things.

\section{A new twist: worst-case switching and the maximum principle}\label{s-margaliot}

{In September of 2003, I received an email from Michael Margaliot, a professor at Tel Aviv University. We had never met or communicated before, although I had seen some of his published work. In his email, Michael said that he had some ideas for solving the above open problem.} The approach that he outlined was based on the concept of \ital{worst-case switching}, an idea that we already encountered in Section~\ref{s-linear-stab} and that had appeared elsewhere in the literature, most notably in the work of Pyatnitskiy and Rapoport~\cite{pyatnitskiy-rapoport-96}. Inspired by that work, Michael proposed to consider an auxiliary control system whose trajectories contain those of the switched system, and to formulate, for this control system, an optimal control problem that consists in driving the state as far away from the origin as possible in a given amount of time. We can then study this optimal control problem with the help of the maximum principle and, under suitable conditions, hope to show that optimal controls are \ital{bang-bang}, i.e., take values only at the extreme points of the control set. Furthermore, we can hope to derive an upper bound on the total number of switches of the optimal controls.
We already know from Sections~\ref{s-commuting} and~\ref{s-gurvits} that having an upper bound on the number of switches---combined with asymptotic stability of the individual modes---leads quite directly to asymptotic stability under switching. In the present case, since optimal controls correspond to the worst-case (most destabilizing) switching, we can conclude as in Section~\ref{s-linear-stab} that the switched system is GUAS.

At the first glance, the approach just described makes no mention of commutation relations, so how is it relevant to our problem? The answer lies in the fact that the optimal control is determined by the sign of a certain function, and switches occur when this function changes sign. As it turns out, the derivatives of this function involve Lie brackets of the vector fields that define the switched system (and the auxiliary control system). When certain Lie brackets vanish, the above function can be shown to be polynomial, and we have a bound on the number of its sign changes. So, there is in fact a direct connection with the Lie brackets.

{Needless to say, Michael's email blew my mind. I had been so focused on solvable Lie algebras and triangular structure that I had completely abandoned an earlier approach, followed by Gurvits, which centered on reachability with a bounded number of switches.
In the linear case, that earlier approach had led me to a dead end while the approach based on triangular structure had proved more fruitful. However, in the nonlinear case the latter approach had led me to a dead end as well. Michael's novel idea was, in essence, to return to the older approach, complemented by the bang-bang principle of optimal control.}

Michael himself had only worked out a particular case of a switched linear system and low-order nilpotency, thus reproving known results. He wrote to me to ask for my help with the nonlinear setting. In order to explain the results that we eventually obtained, it is easiest to begin with a familiar special case of two modes ($m=2$) and second-order nilpotency. This is the direct nonlinear counterpart of the linear condition~\eqref{e-2nd-order-nilpotent-linear} for which Gurvits proved his result in~\cite{gurvits-nilpotent-laa95}, and the simplest next case to consider after the commuting one treated in Section~\ref{s-commuting-nonlinear}. In other words, let us first see how we can prove the following:

\smallskip
\emph{If two globally asymptotically stable vector fields $f_1$ and $f_2$ satisfy
\begin{equation}\label{e-2nd-order-nilpotent-nonlinear}
[f_1,[f_1,f_2]](x)=[f_2,[f_1,f_2]](x)=0 \qquad \foral x\in\R^n
\end{equation}
then the switched nonlinear system~\eqref{e-switched-nonlinear} is GUAS.}
\smallskip

The first step is to define the control system
\begin{equation}\label{e-control-system}
\dot x=f(x)+g(x)u
\end{equation}
with $f:=f_1$, $g:=f_2-f_1$, and the control set $U:=\{0,1\}$. It is clear that for piecewise constant controls $u(\cdot)$ taking values in $U$, the trajectories of~\eqref{e-control-system} coincide with those of the switched system~\eqref{e-switched-nonlinear}. For technical reasons, it is preferable to enlarge the control set to $\bar U:=[0,1]$ and to allow more general control functions $u(\cdot)$ taking values in $\bar U$. The solutions of the resulting control system coincide with those of the differential inclusion
\begin{equation}\label{e-diff-incl}
\dot x\in\text{co}\{f_1(x),f_2(x)\}
\end{equation}
(here `co' denotes the convex hull); they include all solutions of the original switched system.

Now, we pose the following optimal control problem: for an arbitrary given initial condition $x(0)$ and a given time horizon $t_f>0$, find a control $u(\cdot)$ that maximizes the functional
$$
J(u):=|x(t_f)|^2
$$
where $x(\cdot)$ is the state trajectory generated by $u(\cdot)$. Intuitively, we are looking for the worst-case (the most destabilizing) control. If we can show that the resulting
closed-loop system is asymptotically stable, then the same property should hold
for all other controls, and global asymptotic stability of the differential inclusion~\eqref{e-diff-incl}---hence in particular GUAS of the switched
system~\eqref{e-switched-nonlinear}---will be
established.

This problem can be studied with the help of the \ital{maximum principle} of optimal control.\footnote{The reader can find the statement of the maximum principle in almost any textbook on optimal control theory; \cite{cvoc} is my personal favorite.} To this end, we introduce the \ital{Hamiltonian}
$$
H(x,u,p):=p^Tf(x)+p^Tg(x)u
$$
where $p:[0,t_f]\to\R^n$ is the \ital{costate} satisfying the \ital{adjoint equation} $\dot p=-\partial H/\partial x$. The maximum principle then says that
an optimal control must maximize the Hamiltonian pointwise in time. More precisely, at each time $t$, an optimal control should maximize the function $u\mapsto H(x(t),u,p(t))$, where $x(\cdot)$ is the corresponding optimal state trajectory and $p(\cdot)$ is a costate trajectory.

In view of the affine dependence of the right-hand side of the system~\eqref{e-control-system}, and hence of the Hamiltonian $H$, on $u$, it is easy to see how to choose $u$ to maximize $H$. If we define the function
$$
\varphi(t):=p^T(t)g(x(t))
$$
then an optimal control must satisfy
$
u(t)=1$
if $\varphi(t)>0$, and
$u(t)=0$ if $\varphi(t)<0$.

We see that the switches of an optimal control are governed by the sign of the function $\varphi$. Here the announced link with Lie brackets finally appears. Using the definition of $\varphi$ and the differential equations for $x$ and $p$, we can compute the time derivative of $\varphi$ to be
$$
\dot\varphi(t)=p^T(t)[f,g](x(t))
$$
and its second derivative to be
\begin{equation}\label{e-varphi-ddot}
\ddot\varphi(t)=p^T(t)[f,[f,g]](x(t))+p^T(t)[g,[f,g]](x(t))u.
\end{equation}
Recalling~\eqref{e-2nd-order-nilpotent-nonlinear}, it is not hard to see that the second-order Lie brackets in the last expression are 0, forcing $
\ddot\varphi$ to be 0. This of course means that $\varphi$ is a linear function of time, and so it has at most one sign change. Therefore, optimal controls are bang-bang with at most one switch! From this, we can conclude as explained earlier that optimal state trajectories go to 0, and GUAS of the switched system follows.

In the above reasoning, we glossed over one important possibility: what if $\varphi$ is \ital{identically zero}? If this happens, the maximum principle no longer guarantees that an optimal control only takes the values~0 and~1. While we cannot rule out such behavior, fortunately we can use a lemma proved by Hector Sussmann in~\cite{sussmann-bang-bang-siam-1979} to show that in this case, another optimal control exists which is bang-bang but may have one additional switch, i.e, it has at most two switches in total. As we know, with this extra switch the stability proof still goes through.

To go beyond the second-order nilpotency condition~\eqref{e-2nd-order-nilpotent-nonlinear}, let us suppose, referring to~\eqref{e-varphi-ddot}, that we still have $[g,[f,g]](x)=0$ for all $x$, but $[f,[f,g]]$ is nonzero. Then we can proceed to calculate the third time derivative of $\varphi$:
\begin{equation}\label{e-varphi-dddot}
\dddot\varphi(t)=p^T(t)[f,[f,[f,g]]](x(t))+p^T(t)[g,[f,[f,g]]](x(t))u.
\end{equation}
If the third-order Lie brackets in this expression vanish, then $\varphi$ is a quadratic function of time, and so it has at most two sign changes (unless it is identically 0, but we can handle that case in the same way as before). This still gives us a desired upper bound on the number of switches of optimal controls, implying GUAS of the switched system. In a similar fashion, we can formulate Lie-algebraic conditions guaranteeing that $\varphi(\cdot)$ is a polynomial of degree 3, 4, and so on. Unfortunately, these conditions are not quite the $k$th-order nilpotency of the Lie algebra generated by $f$ and $g$, because we need to assume that the lower-order Lie brackets appearing in the $u$-dependent terms in~\eqref{e-varphi-ddot}, \eqref{e-varphi-dddot}, ... are 0 as well.

So far we have only considered the case of $m=2$ modes. For a general $m$, the affine control system corresponding to the switched system~\eqref{e-switched-nonlinear} takes the form
\begin{equation}\label{e-control-affine}
\dot x=f(x)+\sum_{i=1}^{m-1}g_i(x)u_i
\end{equation}
where $f:=f_1$; $g_i:=f_{i+1}-f_1$, $i=1,\dots,m$; $U\in\R^{m-1}$ consists of the standard unit vectors (with one component equal to 1 and the other components equal to 0) and the origin; and $\bar U=\text{co}(U)$ is the standard simplex with these $m$ vertices. Switches of each component $u_i$ of an optimal control are now governed by sign changes of a corresponding function $\varphi_i$, and we can still formulate Lie-algebraic conditions under which these functions are all polynomial. If $d$ is an upper bound on the degrees of these polynomials, then it can be shown that optimal controls have at most $(d+2)^{m-1}-1$ switches. (For $m=2$ and $d=1$ we recover 2 switches as above.)

Full details of this general formulation can be found in our paper with Michael~\cite{margaliot-liberzon-scl}. This approach was also further explored and extended in a follow-up paper that Michael wrote with his MS student Yoav Sharon~\cite{sharon-margaliot}, who subsequently did his PhD with me (on unrelated topics). While these results are encouraging, I still consider the general nonlinear problem to be largely open.

It seems appropriate here to go back to the linear case for a moment and mention another relevant result by Michael.  As we discussed in Section~\ref{s-gurvits}, Gurvits claimed in~\cite{gurvits-nilpotent-laa95} that if $m$ matrices generate a Lie algebra nilpotent of order $k$, then there exists a positive integer $N=N(m,k)$ such that an arbitrary product of exponentials of these matrices can be written as a product of at most $N$ exponential terms. Gurvits himself only proved the case $m=k=2$, in which $N=5$ as expressed by~\eqref{e-gurvits-factorization}. In the paper~\cite{margaliot-counterexample-gurvits}, Michael actually disproved Gurvits' general claim by giving a counterexample where $m=2$, $k=3$, and an $N$ with the above property does not exist. So, my earlier efforts to prove Gurvits' claim were futile. The approach proposed by Michael and developed (for nonlinear systems) in our paper~\cite{margaliot-liberzon-scl} still draws quite heavily on that of Gurvits, since both invoke reachability with a bound on the number of switches. However, the worst-case switching---formalized via optimal control and analyzed using the maximum principle---was a crucial missing ingredient that Michael supplied and that allowed us to make progress.

\section{A throwback}

{In 2009 I met Yuliy Baryshnikov, who a couple of years later became my colleague at University of Illinois. I told him about the research I had been doing earlier on Lie algebras and stability of switched systems, and he surprised me by saying that, yes, he knew that a switched linear system is stable under arbitrary switching if the corresponding matrix Lie algebra is solvable, or ``solvable plus compact" as in Section~\ref{s-compact-levi-factor}. This seemed odd because Yuliy had not really been working in control theory since his days at the Institute of Control Problems in Moscow (where both he and I grew up) in the 1980s, and it was unlikely that he would have read my papers. I asked him where he had seen these results, and he pointed me to two papers by Sergey Kutepov~\cite{kutepov-82,kutepov-84} dating back to 1982 and 1984.}

{It was understandable why I had not come across these references earlier. First, they appeared in a Russian-language journal and were never translated into English.\footnote{The journal was actually Ukrainian. As I'm writing this article, Russia is waging a brutal war in Ukraine, but at the time the distinction hardly mattered.} Of course, I could read Russian perfectly well, but this explains why these papers have not been cited in the Western control-theoretic literature. Second, Kutepov's work did not mention switched systems at all. Instead, it was concerned with exponential stability, uniform over all controls, of bilinear control systems of the form
$
\dot x=\sum_{i=1}^mu_iA_ix
$.
For suitable choices of the control set, such a system captures the behavior of the switched linear system generated by the matrices $A_1,\dots,A_m$; this is a special case of the relationship that we already saw in Section~\ref{s-margaliot} between the affine control system~\eqref{e-control-affine} and the switched nonlinear system~\eqref{e-switched-nonlinear}.}

{These differences notwithstanding, Kutepov's papers~\cite{kutepov-82} and~\cite{kutepov-84} derived the same basic Lie-algebraic stability criteria as the ones in our papers~\cite{lie} and~\cite{agrachev}, respectively. They did not contain constructions of quadratic common Lyapunov functions or the more advanced results of~\cite{agrachev}. On the other hand, Kutepov's proofs were in some places more elegant, revealing a firmer grasp of the theory of Lie algebras than the one I had.}

{There is a joke that I heard from Andy Teel, and so I will attribute it to him. It goes like this: Russian technical literature is not causal; you prove a theorem, and some time later it appears in an old Russian paper. Kutepov's papers from the 1980s, which ``appeared" about 10 years after our work, are a perfect example of this.}

{Of course, in subsequent citations I give full credit to Kutepov. We also submitted a note to Systems and Control Letters explaining the relationship between~\cite{kutepov-82} and~\cite{lie}. For some reason they never published it, but it is posted on my website. Through my Moscow contacts I was eventually able to find Kutepov's email address and wrote to him about how I had been unwittingly following in his footsteps. He was very gracious and thanked me for helping to finally bring proper recognition to his work.}

{Had I learned about Kutepov's work soon after writing my papers containing the same results, I would have certainly been devastated. In that sense, I am thankful that many years had passed, I had a faculty job, and my research reputation was less dependent on that early work. I was also mature enough to understand that our collective knowledge is more important than who did what first.}

\section{Robustness of Lie-algebraic stability conditions}\label{s-robustness}

Lie-algebraic stability conditions, while mathematically appealing and often fairly easily checkable, suffer from one serious drawback: they are not robust with respect to small perturbations of the system data. For example, if we take two matrices that commute with each other and perturb one of them slightly, they will cease to commute. If we take a family of matrices generating a nilpotent or solvable Lie algebra and introduce arbitrarily small errors in their entries, the new Lie algebra will no longer possess any helpful structure (see~\cite[Section~A.6]{agrachev} for a precise result along these lines).

On the other hand, the stability properties that these Lie-algebraic conditions establish \ital{are} robust to small perturbations. For the switched linear system~\eqref{e-switched-linear}, this is especially easy to see with the help of a quadratic common Lyapunov function characterized by the inequalities~\eqref{e-qclf-PQ}, which exists in all situations discussed in Sections~\ref{s-commuting}--\ref{s-compact-levi-factor}. Suppose that we are given a family of perturbed matrices of the form $\bar A_p=A_p+\Delta_p$, $\pp$. Then $V(x)=x^TPx$ is still a common Lyapunov function for the perturbed systems $\dot x=\bar A_px$ if the perturbations matrices $\Delta_p$ satisfy
\begin{equation}\label{e-robustness-bound}
\|\Delta_p\|_2<\frac{\lambda_{\min}(Q)}{2\lambda_{\max}(P)}
\end{equation}
where $\|\cdot\|_2$ is the matrix norm induced by the Euclidean norm and $\lambda_{\min}(\cdot)$ and $\lambda_{\max}(\cdot)$ denote the smallest and the largest eigenvalue of a symmetric matrix, respectively (see~\cite[p.\ 342]{khalil-book-3ed} or~\cite[p.\ 42]{book}). This suggests that, instead of requiring that suitable commutators vanish, robust Lie-algebraic stability conditions should ask that these commutators be sufficiently small, in order to guarantee that a given family of matrices is close to a family of commuting matrices or to a family of matrices generating a nilpotent, solvable, or ``solvable plus compact" Lie algebra.

{I already explained one consequence of my conversation with Yuliy Baryshnikov in 2009. Another, perhaps more important outcome of our discussions was that Yuliy suggested several approaches for developing robust Lie-algebraic stability conditions for switched linear systems.} One approach, carried out in discrete time, involves rearranging the order of matrices in a product (somewhat in the spirit of Gurvits' ``elementary shufflings" mentioned in Section~\ref{s-gurvits}) and characterizing
the effect of such a rearrangement in terms of the norms of the commutators of the matrices. The switched linear system is GUES if the commutators are small enough, which gives a robust version of the result from Section~\ref{s-commuting}. Another approach is in continuous time and utilizes the structure of the Lie algebra. In the notation of Section~\ref{s-compact-levi-factor}, we showed GUES when the matrices in the subalgebra ${\mathfrak{s}}$
have sufficiently small norms compared to the stability margin of
the matrices in the solvable part ${\mathfrak{r}}$, which is
 a robust version of the
result from Section~\ref{s-solvable}. We also showed that GUES is preserved
when noncompact perturbations are introduced, as long as they are
small compared to the real parts of the eigenvalues of the
matrices in the solvable part ${\mathfrak{r}}$; this is a robust version of the stability criterion from Section~\ref{s-compact-levi-factor}. All of these results are documented in our 2012 paper~\cite{ABL-SCL} with Yuliy and Andrei Agrachev.
In the follow-up paper~\cite{Loja}, Yuliy and I used an inequality due to  {\L}ojasiewicz to relate the size of commutators of given matrices to the distance from this family of matrices to a family of commuting matrices (or matrices generating a nilpotent or solvable Lie algebra). Coupled with the bound~\eqref{e-robustness-bound}, this again leads to robust stability criteria involving bounds on the commutators.

{The discrete-time approach from~\cite{ABL-SCL} was later taken up by Atreyee Kundu, who did her PhD under Debasish Chatterjee, my former PhD student. In the paper~\cite{atreyee-kundu-scl-20} Atreyee adapted Yuliy's method to switching signals satisfying a dwell-time condition, while in~\cite{atreyee-debasish-automatica-20} Atreyee and Debasish further extended the results by allowing the presence of unstable modes.} Another relatively recent work~\cite{gil-small-commutators-2019} derives an upper bound on the commutator of two Hurwitz matrices under which the procedure~\eqref{e-nb94}--\eqref{e-VPm}, with $m=2$, still yields a quadratic common Lyapunov function.

At the end of Section~\ref{s-solvable} we mentioned the paper~\cite{haimovich-braslavsky-tac11} which addressed stabilization under arbitrary switching via finding controller gains that achieve simultaneous triangularization of the closed-loop matrices. The same authors---Haimovich and Braslavsky---also developed a robust version of their procedure, which only requires approximate simultaneous triangularization~\cite{haimovich-braslavsky-cdc10}.

\section{Confusions}

%

Reading the above story, one might get the impression of clear vision conceived, systematic efforts expended, and steady progress achieved. For me, this could not be farther from the truth. The part of the work described here in which I was personally involved was in reality a fairly random sequence of confused attempts, unexpected turns, and frustrating failures. I was also working on several other research topics, quite unrelated to the problem addressed in this article, and there were long periods when I was not really thinking about this problem until some spontaneous conversation would nudge me to return to it. It is only \ital{after} the work has been done that a coherent story has gradually emerged. I have been considering for some time the possibility of writing a technical survey on this subject, but it was only very recently---as I was teaching my course on switched systems and telling various stories to the students---that the idea occurred to me to write this article from an informal personal perspective. I hope that some readers will find this account useful, and the missing technical details can always be found in the cited papers. I should also stress that, since this article was conceived as a personal story and not as an exhaustive survey of relevant research, my overview of the literature is far from being complete.


\end{document}